\documentclass[11pt]{amsart}

\usepackage{amsmath, amssymb, amsfonts, amsthm, amscd}
\usepackage{manfnt}
\usepackage{mathtools}
\usepackage{fullpage}
\usepackage[all]{xy}
   \SelectTips{cm}{10}
\usepackage{txfonts}

\newtheorem{thm}{Theorem}[section]
\newtheorem{lemma}[thm]{Lemma}
\newtheorem{prop}[thm]{Proposition}
\newtheorem{cor}[thm]{Corollary}

\newtheorem{prop-conj}[thm]{Proposition-Conjecture}

\theoremstyle{definition}
\newtheorem{defn}[thm]{Definition}

\theoremstyle{remark}
\newtheorem{rmk}[thm]{Remark}

\theoremstyle{remark}

\theoremstyle{remark}

\theoremstyle{remark}
\newtheorem{eg}[thm]{Example}

\newcommand{\Q}{\mathbb{Q}}

\newcommand{\Z}{\mathbb{Z}}
\newcommand{\CC}{\mathbb{C}}
\newcommand{\RR}{\mathbb{R}}

\newcommand{\Qlb}{\overline{\mathbb{Q}}_\ell}
\newcommand{\af}{\mathbf{A}_F}
\newcommand{\afp}{\mathbf{A}_{F}}

\DeclareMathOperator{\Hom}{Hom}
\DeclareMathOperator{\Ext}{Ext}

\DeclareMathOperator{\BC}{BC}
\DeclareMathOperator{\Sym}{Sym}

\DeclareMathOperator{\rec}{rec}

\DeclareMathOperator{\coker}{coker}

\DeclareMathOperator{\diag}{diag}

\newcommand{\gal}[1]{\Gamma_{#1}} 
\newcommand{\Gal}{\mathrm{Gal}} 

\newcommand{\into}{\hookrightarrow}
\newcommand{\onto}{\twoheadrightarrow}

\newcommand{\mc}{\mathcal}

\newcommand{\mr}{\mathrm}
\newcommand{\mbf}{\mathbf}
\newcommand{\mbb}{\mathbb}

\newcommand{\tG}{\widetilde{G}}
\newcommand{\tT}{\widetilde{T}}
\newcommand{\tZ}{\widetilde{Z}}
\newcommand{\tH}{\widetilde{H}}
\newcommand{\tpi}{\tilde{\pi}}
\newcommand{\tomega}{\tilde{\omega}}
\newcommand{\tmu}{\tilde{\mu}}

\newcommand{\LtG}{{}^L \widetilde{G}}

\author{Stefan Patrikis}
\address{Department of Mathematics \\  MIT \\ Cambridge, MA 02139}
\email{patrikis@math.mit.edu}
\title{On the sign of regular algebraic polarizable automorphic representations}
\date{July 2013}
\thanks{This paper owes its existence to questions raised by Richard Taylor and Claus Sorensen: Richard asked the question addressed in Theorem \ref{mainthm}, and Claus stressed that Corollary \ref{pm1liftstotreal} should be true. I am very grateful to both of them. This work was carried out while a member of the Institute for Advanced Study, supported by NSF grant DMS-1062759.}
\subjclass[2010]{11F80,11R39.}
\keywords{geometric Galois representation, algebraic automorphic representations}

\begin{document}
\begin{abstract}
 We remove a parity condition from the construction of automorphic Galois representations carried out in the Paris Book Project. We subsequently generalize this construction to the case of `mixed-parity' (but still regular essentially self-dual) automorphic representations over totally real fields, finding associated geometric projective representations. Finally, we optimize some of our previous results on finding geometric lifts, through central torus quotients, of geometric Galois representations, and apply them to the previous mixed-parity setting. 
\end{abstract}
\maketitle

\section{Introduction}\label{sectionintro}
This note discusses a few closely-related parity phenomena arising in the study of automorphic Galois representations over totally real fields. The prototype of these phenomena is Weil's distinction, in his seminal investigation \cite{weil:characters} of `algebraic' Hecke characters, between type $A_0$ and type $A$ Hecke characters of number fields. Both have (incomplete) L-series with algebraic coefficients, but only the former give rise to compatible systems of $\ell$-adic Galois characters; for this reason the latter, and their analogues in higher rank, have been somewhat neglected. But keeping them in mind provides useful intuition for a number of natural arithmetic questions; here we briefly summarize the three addressed in this note:
\begin{itemize}
\item We remove a parity condition from the construction of automorphic Galois representations for RAESDC automorphic representations over totally real fields (due to many people; for a review, see \cite{blggt:potaut}; for our result, see Theorem \ref{mainthm} below). Whereas these Galois representations are constructed via a descent to appropriate unitary groups, to understand the parity condition we make use of a descent to quasi-split $\mr{GSpin}$ groups, using work of Asgari-Shahidi and Hundley-Sayag (see \cite{asgari-shahidi:splitGSpin}, \cite{asgari-shahidi:nonsplitGSpin}), Hundley-Sayag (\cite{hundley-sayag:GSpindescent}, and \cite{hundley-sayag:GSpindescentproofs}).
\item We extend the above-mentioned construction of automorphic Galois representations from the L-algebraic to the `mixed-parity' case, finding associated geometric \textit{projective} representations; the same technique yields such Galois representations associated to certain automorphic representations of $\mr{Spin}$ groups. (See Theorem \ref{mixedgalois} and Corollary \ref{pm1liftstotreal} for these `mixed-parity' cases; as a much easier warm-up, and some evidence for the Buzzard-Gee conjecture, see Proposition \ref{galoisreps} for the L-algebraic case on $\mr{GSpin}$ groups.) To handle the mixed-parity case, a new idea is required to establish a `projective' variant of the patching lemma of Blasius-Rogawski (\cite{blasius-rogawski:motiveshmfs}; see too \cite{sorensen:patching}). 
\item Prompted by a question of Claus Sorensen, we show that these projective representations have geometric lifts to $\mr{GL}_N(\Qlb)/\{\pm 1\}$-valued representations (but not to geometric $\mr{GL}_N(\Qlb)$-valued representations!). To do this we streamline and optimize some of the general Galois-theoretic lifting results of \cite[\S 3.2]{stp:variationsarxiv}. (See Proposition \ref{reallift} and Corollary \ref{pm1liftstotreal}.)
\end{itemize}
See the individual sections for more context and explanation of these problems; for more general background on these sorts of parity questions, see \cite{stp:variationsarxiv}.

Before continuing to the main body of the paper, we also review a little terminology. For a number field $F$ and $G/F$ a connected reductive group, we have the following notions of algebraicity for an automorphic representation $\pi$ of $G(\af)$. For all $v \vert \infty$, the archimedean local Langlands correspondence yields a representation
\[
\phi_v= \rec_v(\pi_v) \colon W_{F_v} \to {}^L G
\]
of the Weil group $W_{F_v}$ to an $L$-group for $G$. Fixing a maximal torus $T^\vee$ of the dual group $G^\vee \subset {}^L G$, and choosing an algebraic closure $\overline{F}_v$ and an isomorphism $\iota_v \colon \overline{F}_v \xrightarrow{\sim} \CC$, we may assume (up to $G^\vee$-conjugation) that $\phi_v|_{W_{\overline{F}_v}}$ has the form
\[
z \mapsto \iota_v(z)^{\mu_{\iota_v}} \bar{\iota}_v(z)^{\nu_{\iota_v}} \in T^\vee,
\]
for ($\CC$-linear combinations of) co-characters $\mu_{\iota_v}, \nu_{\iota_v} \in X_{\bullet}(T^\vee) \otimes_{\Z} \CC$ satisfying $\mu_{\iota_v}- \nu_{\iota_v} \in X_{\bullet}(T^\vee)$. The best hope (see for instance \cite{buzzard-gee:alg}) is that $\pi$ will have associated ${}^L G(\Qlb)$-valued Galois representations when the various $\mu_{\iota_v}$ and $\nu_{\iota_v}$ (for all $v \vert \infty$) all in fact belong to $X_{\bullet}(T^\vee)$; intuitively, these are those for which we can `see a Hodge structure' on the hoped-for motive. Following Buzzard and Gee, we call $\pi$ satisfying this archimedean condition `L-algebraic.' Although we do not make use of it, we mention also the `C-algebraic' condition-- terminology again due to Buzzard and Gee-- which describes those $\pi$ for which the $\mu$'s and $\nu$'s all lie in $\rho+X_{\bullet}(T^{\vee})$, for $\rho$ the usual half-sum of positive roots. Both L- and C-algebraicity generalize Weil's notion of type $A_0$ Hecke character (for any $G= \mr{GL}_{2n+1}$, including $\mr{GL}_1$, they give the same notion), although in this paper we are equally concerned with the more general notion of type $A$ Hecke character, and its higher-rank analogues. Here again there are various notions one might use. The broadest-- and the direct analogue of Weil's type $A$ condition-- is to require that the $\mu$'s and $\nu$'s simply be quasi-co-characters, i.e. lie in $X_{\bullet}(T^\vee)_{\Q}$. One might expect that this condition prescribes exactly those automorphic representations with algebraic Satake parameters. But there are fundamental constraints, stemming from the Ramanujan conjecture, on the possible infinity-types of sufficiently `non-degenerate' automorphic representations, and it is often useful to build some of these constraints into a definition. This motivates the notion `W-algebraic,' describing those $\pi$ for which all the $\mu$'s and $\nu$'s lie in $\frac{1}{2}X_{\bullet}(T^\vee)$; see \cite{stp:variationsarxiv} for more discussion of these matters. 

Finally, in this paper, at least on the group $\mr{GL}_N$, we will want to draw attention to a certain subset of those $\pi$ that are $W$-algebraic but not $L$-algebraic (nor twists of $L$-algebraic representations):
\begin{defn}\label{mixedparitydef} 
An automorphic representation $\pi$ of $\mr{GL}_N(\af)$ will be said to be `mixed-parity' if, for some proper subset of $\{v \vert \infty\}$, the $\mu_{\iota_v}$ and $\nu_{\iota_v}$ lie in $X_{\bullet}(T^\vee)$ (in which case we will also say $\pi_v$ is L-algebraic); while for the complementary, still proper, subset of $\{v \vert \infty\}$, $\pi_v \otimes |\cdot|^{1/2}$ is L-algebraic.
\end{defn} 
This definition simply extends the classical terminology of `mixed-parity' Hilbert modular forms. See \S \ref{sectionmixed} for more discussion of this notion: it is only relevant for $N$ even (see Lemma \ref{nomixedodd}), in which case it can-- perhaps tellingly-- be rephrased as: $\pi_v$ is L-algebraic for some $v \vert \infty$, and C-algebraic for other $v \vert \infty$.
\section{Removing a sign condition from the construction of certain automorphic Galois representations}
Let $F$ be a totally real number field, and let $\Pi$ be a regular, (C- or L-) algebraic, essentially self-dual cuspidal automorphic representation of $\mr{GL}_N(\afp)$, so for some (type $A_0$) Hecke character $\omega$ we have $\Pi \cong \Pi^{\vee} \otimes \omega$. One knows (by the work of many people; see \cite[Theorem 2.1.1]{blggt:potaut} for a resume) how to associate automorphic Galois representations to such $\Pi$ under the additional hypothesis that the sign $\omega_v(-1)$ is independent of $v \vert \infty$ (in the terminology \cite{blggt:potaut}, such $\Pi$ are `polarizable'). Our first result shows that this additional hypothesis is in fact superfluous:
\begin{thm}\label{mainthm}
Let $F$ be a totally real field, and let $\Pi$ be a regular, C-algebraic or L-algebraic, cuspidal automorphic representation of $\mr{GL}_N(\af)$ satisfying a self-duality $\Pi \cong \Pi^\vee \otimes \omega$ for some Hecke character $\omega \colon C_F= \mathbf{A}_F^\times/F^\times \to \CC^\times$. Then $\omega_v(-1)$ is independent of $v \vert \infty$, and consequently one can associate compatible systems of $\ell$-adic representations to $\Pi$, as in \cite[Theorem 2.1.1]{blggt:potaut}.
\end{thm}
The result is obvious when $N$ is odd, so from now on we let $N=2n$ be even. The key ingredient in the proof of theorem is the descent of $\Pi$ to a suitable $\mr{GSpin}$ group, thanks to the work of Asgari-Shahidi (\cite{asgari-shahidi:splitGSpin} and \cite{asgari-shahidi:nonsplitGSpin}) and Hundley-Sayag (\cite{hundley-sayag:GSpindescent} and \cite{hundley-sayag:GSpindescentproofs}); there the sign $\omega_v(-1)$ of interest can be interpreted as a central character, where it is determined by the `parity' of the corresponding discrete series representation at $v$.
\begin{eg}\label{mixedparityeg}
The general argument is modeled on the case $n=1$ (note that $\mr{GSpin}_3 \cong \mr{GL}_2$), where it is nearly a tautology (so much so that the relevant generalization may be obscured). In this case, $\Pi$ is a Hilbert modular representation with central character $\omega$, and if we let $k_v$ be its `weight' at $v \vert \infty$ in the sense of classical modular forms, then (having normalized $\Pi$ to be unitary) the archimedean L-parameters $\phi_v \colon W_{\RR} \to \mr{GL}_{2}(\CC)$ have the following form
\begin{align*}
&z \mapsto 
\begin{pmatrix}
(z/\bar{z})^{\frac{k_v-1}{2}} & 0 \\
0 & (z/\bar{z})^{\frac{1-k_v}{2}}
\end{pmatrix} \\
&j \mapsto
\begin{pmatrix}
0 & (-1)^{k_v-1} \\
1 & 0
\end{pmatrix}.
\end{align*}
Thus $\omega_v(-1)= \det \phi_v(j)= (-1)^{k_v}$. The crucial point is that for $\Pi$ to be C-algebraic (respectively, L-algebraic), all $k_v$, for $v \vert \infty$, must be even (respectively, odd). Thus, under the assumptions of the theorem, $\omega_v(-1)$ is independent of $v \vert \infty$. Those Hilbert modular representations $\Pi$ for which $\omega_v(-1)$ varies with $v$, the so-called `mixed-parity' representations, are $W$-algebraic in the sense of \cite{stp:variationsarxiv}, and should be thought of as higher rank analogues of type $A$ but not $A_0$ Hecke characters in the sense of Weil (\cite{weil:characters}).
\end{eg}


We now recall the deep results on generic transfer and descent for $\mr{GSpin}$ groups that allow this simple argument to be extended to higher rank. Let $\tG$ denote a quasi-split general spin group over $F$; later on we will reserve $G$ for the corresponding $\mr{Spin}$ group. We choose a based root datum, along with a splitting, and form the associated L-group $\LtG$.  The details of these choices will not be too important for us, so we refer the reader to \cite[\S 2.1]{asgari-shahidi:nonsplitGSpin} for explicit descriptions. In this section we will carry out explicit matrix rather than root-theoretic calculations, so to be clear: our convention here and throughout is that $\mr{GSO}_{2n}(\CC)$ will be defined with respect to the symmetric pairing
$\begin{pmatrix}
 & 1_n \\
1_n & 
\end{pmatrix}$, while $\mr{GSp}_{2n}(\CC)$ will be defined with respect to the alternating pairing 
$\begin{pmatrix}
 & 1_n \\
-1_n & 
\end{pmatrix}$. In all cases, then, the dual group $\tG^\vee \subset \mr{GL}_{2n}(\CC)$ has a diagonal maximal torus of the form 
\[
\tT^\vee= \{ \diag(t_1, \ldots, t_n, t_1^{-1} x, \ldots, t_n^{-1} x): x, t_i \in \CC^\times \}.
\]
$\tG$ is one of the following three types, which we list along with (the Galois form of) its $L$-group ${}^L \tG$:
\begin{itemize}
\item the split group $\mr{GSpin}_{2n+1}/F$, with $\LtG= \mr{GSp}_{2n}(\CC) \times \gal{F}$;
\item the split group $\mr{GSpin}_{2n}/F$, with $\LtG= \mr{GSO}_{2n}(\CC) \times \gal{F}$;
\item one of the quasi-split but not split groups $\mr{GSpin}_{2n}^{\mu}/F$ associated to a quadratic extension $F'/F$ cut out by a character $\mu \colon C_F= \mathbf{A}_F^\times/ F^{\times} \to \{ \pm 1 \}$. We can choose the based root datum and splitting so that $\LtG= \mr{GSO}_{2n}(\CC) \rtimes \gal{F}$, with the action of $\gal{F}$ factoring through $\Gal(F'/F)$, where it is given by conjugation by the matrix
\[
h= \left(
\begin{array}{cc|cc}
1_{n-1} & & & \\
& 0 &  & 1 \\ \hline
&  & 1_{n-1} & \\
& 1 & & 0
\end{array} \right).
\]
(See \cite[\S 2.2]{asgari-shahidi:nonsplitGSpin}; note that their $h$ looks a little different from ours since they take an orthogonal pairing for which a diagonal maximal torus in $\tG^\vee$ has the form $\diag(t_1, \ldots, t_n, t_n^{-1}x, \ldots, t_1^{-1}x)\}$.) Note that $h$ lies in $\mr{GO}_{2n}(\CC)$ but not in $\mr{GSO}_{2n}(\CC)$; $h$-conjugation preserves our (implicit) choice of based root datum, swapping two simple roots.
\end{itemize}

The transfer of interest is with respect to the $L$-homomorphism $\iota_{\tG} \colon \LtG \to \mr{GL}_{2n}(\CC) \times \gal{F}$ (the target of course being the $L$-group of $\mr{GL}_{2n}/F$) given as follows: in the split cases, $\iota_{\tG}$ is the obvious extension to $\LtG$ of the standard $2n$-dimensional representation of $\tG^{\vee}$, and in the non-split case, $\iota_{\tG}$ is given by:\footnote{Note that there is a typo in the definition of $\iota_{\tG}$ in \cite[3.3]{asgari-shahidi:nonsplitGSpin}; they refer to \cite{cogdell-ps-shahidi:quasisplit}, where the correct definition is given.}
\begin{align*}
\iota_{\tG} \colon \mr{GSO}_{2n}(\CC) \rtimes \gal{F}  \to \mr{GL}_{2n}(\CC) \times \gal{F} \\
(g, \gamma)  \mapsto 
\begin{cases}
(g, \gamma) & \text{if $\gamma|_{F'}= 1$,}\\
(gh, \gamma) & \text{if $\gamma|_{F'} \neq 1$.}
\end{cases}
\end{align*}
In the next theorem, we summarize what we will need from the works of Asgari-Shahidi and Hundley-Sayag; note that the results in those papers are in fact stronger, giving an $L$-function criterion for describing the descent, and extending the descent to isobaric, not merely cuspidal, representations.
\begin{thm}[Asgari-Shahidi, Hundley-Sayag. See Theorem 4.26 of \cite{asgari-shahidi:nonsplitGSpin}]\label{descent}
Let $k$ be any number field, and let $\tpi$ be a unitary\footnote{The unitary assumption in this theorem can easily be removed by twisting.} globally generic cuspidal automorphic representation of $\tG(\mbf{A}_k)$ with central character $\omega= \omega_{\tpi}$. Then $\tpi$ has a unique transfer to an automorphic representation $\Pi$ of $\mr{GL}_{2n}(\mbf{A}_k)$ such that for all infinite places and all finite places $v$ of $k$ at which $\tpi_v$ is unramified and $F'/F$ is unramified, the local L-parameter $\Phi_v$ of $\Pi_v$ is the transfer $\iota_{\tG} \circ \phi_v$ of the local L-parameter $\phi_v \colon W_{k_v} \to \LtG$ of $\tpi_v$. This transfer satisfies
\[
\Pi \cong \Pi^\vee \otimes \omega,
\]
and, writing $\omega_{\Pi}$ for the central character of $\Pi$, we have
\[
\omega_{\Pi}= \omega^n \mu,
\]
where $\mu$ is trivial in the split cases, and is the quadratic character associated to $\tG$ in the quasi-split, non-split cases.\footnote{To be precise, when we regard $\omega$ as a character of $\mbb{A}_k$, we mean the restriction of $\omega$ to $Z_{\tG}^0(\mbb{A}_k)$-- here $Z_{\tG}^0$ denotes the connected component of the center of $\tG$, which is disconnected in the $D_n$ cases. See below.}

Conversely, any unitary cuspidal automorphic representation $\Pi$ of $\mr{GL}_{2n}(\mbf{A}_k)$ satisfying $\Pi \cong \Pi^{\vee} \otimes \omega$ for some Hecke character $\omega$ is the functorial transfer, from one of the above groups $\tG$, of such a $\tpi$.
\end{thm}
Note that the established transfer is compatible with archimedean L-parameters; we will make crucial use of this. Which group $\tG$ provides the automorphic representation $\tpi$ can be read off from properties of the pair $(\Pi, \omega)$: namely, exactly one of the incomplete (throwing away a finite set $S$ of places of $k$ containing the archimedean places) L-functions $L^S(\wedge^2 \otimes \omega^{-1}, \Pi, s)$, $L^S(\Sym^2 \otimes \omega^{-1}, \Pi, s)$, has a pole at $s=1$; in the former case $\Pi$ descends to $\mr{GSpin}_{2n+1}$, and in the latter it descends to $\mr{GSpin}_{2n}$ (when $\omega_{\Pi}/\omega^n=1$) or a quasi-split form (when $\omega_{\Pi}/\omega^n$ is a non-trivial quadratic character). All that matters for us is that $\Pi$ descends to at least one of these groups.

Before giving the proof of Theorem \ref{mainthm}, we recall how central characters of automorphic representations can be computed from local L-parameters, and we make this explicit in the cases of interest.\footnote{This description of the central character in terms of the L-parameter is something that must be proven along with any given case of local Langlands. It is known in general at archimedean and unramified places. See \cite[pg 21-25]{langlands:archllc} and \cite[\S 10]{borel:L}.} For any number field $k$, and any connected reductive group $G/k$, with center $Z_G$ and connected center $Z_G^0$, there is a canonical surjection of $L$-groups ${}^L G \to {}^L Z_G^0$. We can then compose local parameters $W_{k_v} \xrightarrow{\phi_v} {}^L G \to {}^L Z^0_G$, and by the local Langlands correspondence for tori obtain a character $Z^0_G(k_v) \to \CC^\times$. If the center of $G(k_v)$ is contained in $Z_G^0(k_v)$, then this suffices to define the central character. Even when $Z_G$ is disconnected, if we are only interested in the restriction of the central character to $Z_G^0(k_v)$, this description suffices. In the general case, we embed $Z_G$ into a $k$-torus $Z'$ and enlarge $G$ correspondingly, letting
\[
G'= (G \times Z')/Z_G,
\]
so that $Z'$ can be identified with the center of $G'$. Then we \textit{lift}\footnote{That this is always possible, and indeed also for the global Weil group $W_k$, is a theorem of Labesse (\cite{labesse:lifting}); it is essentially an elaboration on Tate's theorem that $H^2(\gal{k}, \Q/\Z)=0$.}  $\phi_v \colon W_{k_v} \to {}^L G$ across the quotient ${}^L G' \to {}^L G$and proceed as before to define a character $Z'(k_v) \to \CC^\times$ whose restriction to $Z_G(k_v)$ is independent of the choice of lift and gives the desired central character.

For us, $G$ will be one of the quasi-split general spin groups that we have been denoting $\tG$. Let $v$ be an archimedean (real) place of $\tG$. Let $j$ denote an element of $W_{F_v}- W_{\overline{F}_v}$ satisfying $j^2=-1$ and $jzj^{-1}= \bar{z}$ for all $z \in W_{\overline{F}_v}$. We then have the explicit descriptions:
\begin{enumerate}
\item $\tG= \mr{GSpin}_{2n+1}/F$. The center $Z_{\tG}$ is a split torus, and the $L$-homomorphism $\LtG \to {}^L Z_{\tG}$ can be identified with the symplectic similitude character $\mr{GSp}_{2n}(\CC) \xrightarrow{\nu} \CC^\times$. Thus, if for $v \vert \infty$ $\tpi_v$ has L-parameter $\phi_v \colon W_{F_v} \to \LtG$, then $\omega_{\tpi_v}(-1)= \nu \circ \phi_v(j)$.
\item $\tG= \mr{GSpin}_{2n}/F$. The center $Z_{\tG}$ is not connected, but $\LtG \to {}^L Z_{\tG}^0$ can be identified with the orthogonal similitude character $\mr{GSO}_{2n}(\CC) \xrightarrow{\nu} \CC^\times$. We still have $\omega_{\tpi_v}(-1)= \nu \circ \phi_v(j)$, where $-1$ refers to $-1 \in Z_{\tG}^0(F_v) \cong \RR^\times$.
\item $\tG= \mr{GSpin}_{2n}^{\mu}/F$. As in the split case: $Z_{\tG}^0$ is still a split torus.
\end{enumerate}

\section{Proof of Theorem \ref{mainthm}}\label{Stwo}
We now prove Theorem \ref{mainthm}. Let $\Pi$ be our given RAESDC automorphic representation; twisting, we may assume that $\Pi$ is unitary; this may replace our C-algebraic $\Pi$ with an L-algebraic $\Pi$, but no matter. We therefore have $\Pi \cong \Pi^\vee \otimes \omega$ for some finite-order Hecke character $\omega$ (any unitary type $A$ Hecke character of a totally real field is finite order). For $v \vert \infty$, the restriction to $W_{\overline{F}_v}$ of the archimedean L-parameter $\Phi_v$ of $\Pi_v$ takes the form\footnote{A choice of isomorphism $\iota_v \colon \overline{F}_v \xrightarrow{\sim} \CC$ will be implicit.}
\[
\Phi_v(z)= \mr{diag} \left( z^{p_{v, 1}} \bar{z}^{q_{v, 1}}, \ldots, z^{p_{v, 2n}} \bar{z}^{q_{v, 2n}} \right) \in \mr{GL}_{2n}(\CC),
\]
where by regularity the various $p_{v, i}$ (for fixed $v$) are distinct, and likewise for the $q_{v, i}$. In fact, $\Pi_{\infty}$ is automatically tempered (see \cite[Lemme 4.9]{clozel:alg}), so (recall we have assumed $\Pi$ is unitary) $p_{v, i}+q_{v, i}=0$ for all $v$ and $i$; invoking the self-duality of $\Pi$ and replacing $\Phi_v$ by a suitable conjugate, we may therefore assume
\[
\Phi_v(z)= \mr{diag} \left( (z/\bar{z})^{p_{v, 1}}, \ldots, (z/\bar{z})^{p_{v, n}}, (z/\bar{z})^{-p_{v, 1}}, \ldots, (z/\bar{z})^{-p_{v, n}} \right),
\]
where each $p_{v, i}$ is non-zero, and $p_{v, i} \neq \pm p_{v, j}$ for any $j \neq i$.

By Theorem \ref{descent}, $\Pi$ descends to an automorphic representation $\tpi$ of one of the $\mr{GSpin}$ groups $\tG$. As before, denote by $\phi_v \colon W_{F_v} \to \LtG$ the L-parameter (defined up to $\tG^\vee$-conjugation) of $\tpi_v$ for all $v \vert \infty$. Up to $\tG^\vee$-conjugation, we may assume $\phi_v|_{W_{\overline{F}_v}}$ lands in $\tT^\vee$. By regularity, the only ambiguity in conjugating $\phi_v$ into $\tT^\vee$ comes from the Weyl group of $\tG^\vee$,
which in all cases does no more than permute the $\{p_{v, i}\}$ and exchange some $p_{v, i}$ with $-p_{v, i}$. Thus we may continue to assume that $\phi_v$ (rather than $\Phi_v$) takes the form
\[
\phi_v(z)= \mr{diag} \left( (z/\bar{z})^{p_{v, 1}}, \ldots, (z/\bar{z})^{p_{v, n}}, (z/\bar{z})^{-p_{v, 1}}, \ldots, (z/\bar{z})^{-p_{v, n}} \right),
\]
although note that these are not necessarily the same $p_{v, i}$ as before-- namely, in the $D_n$ case there can be distinct (up to $\tG^\vee$-conjugation) parameters $\phi_v$ that nevertheless yield the same parameter up to $\mr{GL}_{2n}(\CC)$-conjugation. 

Given such a $\phi_v|_{W_{\overline{F}_v}}$, we now check which extensions to a full L-parameter $\phi_v \colon W_{F_v} \to \LtG$ are possible, and use this to compute $\omega_v(-1)$ in each case. Note that for our given $\Pi$, either all $p_{v, i}$ lie in $\Z$, or all lie in $\frac{1}{2}+ \Z$; the calculation of $\omega_v(-1)$ will turn out only to depend on the group $\tG$ and this integer/half-integer alternative. Let $J \in \tG^\vee$ have the property that $\phi_v(j)=(J, c) \in \LtG$ extends $\phi_v|_{W_{\overline{F}_v}}$ to a well-defined L-parameter. 
\begin{lemma}\label{cent}
The centralizer $\mr{Cent}_{\tG^\vee}(\phi_v(W_{\overline{F}_v}))$ is equal to $\tT^\vee$. In particular, any other $J' \in \tG^\vee$ such that $\phi_v(j)= (J', c)$ gives a well-defined extension of the L-parameter is of the form $J'= Jt$ for some $t \in \tT^\vee$.
\end{lemma}
\proof
A simple calculation, whose details we omit.
\endproof
The next lemma completes the proof of Theorem \ref{mainthm} in the split case:
\begin{lemma}\label{splitcase}
Suppose $\tG$ is one of the split groups $\mr{GSpin}_{2n+1}/F$ or $\mr{GSpin}_{2n}/F$. Having specified $\tG$, the value $\omega_v(-1)$ is determined by whether the $p_{v, i}$ are integers or lie in $\frac{1}{2} +\Z$. In particular, C-algebraicity (or L-algebraicity) of our original $\Pi$ implies that $\omega_v(-1)$ is independent of $v \vert \infty$.
\end{lemma}
\proof
There are four cases to deal with, depending on the group $\tG$, and on whether the $p_{v, i}$ are all integers or all in $\frac{1}{2}+\Z$. In each case we write down a candidate for $J$ and compute the central character from the L-parameter, showing that the value $\omega_v(-1)$ is independent of the choice of $J$. We label the cases $(m, \star)$ where $m = 2n $ or $2n+1$, depending on $\tG$, and $\star= \Z$ or $\frac{1}{2} +\Z$, depending on where the $p_{v, i}$ live.
\begin{enumerate}
\item $(2n+1, \Z)$. Here $J^2= \phi_v(-1)= 1$, so we may take $J= 
\begin{pmatrix}
 & 1_n \\
1_n & 
\end{pmatrix}$. Denote by $\nu \colon \tG^\vee= \mr{GSp}_{2n}(\CC) \to \mbb{G}_m$ the symplectic multiplier. Then $\nu(J)= -1$. Consider another candidate $J' =J t$, with $t \in \tT^\vee$ of the form $\begin{pmatrix}
s & \\
 & s^{-1}x
\end{pmatrix}$, where $s$ is an $n \times n$ diagonal matrix and $\nu(t)=x$. $(J')^2=1$ forces $x=1$, and then we see $\nu(J')=-1$ as well. This forces $\omega_v(-1)= \nu \circ \phi_v(j)= -1$.
\item $(2n+1, \frac{1}{2}+\Z)$. The calculation is similar, except now we have $J^2= \phi_v(-1)= -1$, so we may take $J=
\begin{pmatrix}
 & 1_n \\
 -1_n & 
\end{pmatrix}$. We conclude that $\omega_v(-1)= 1$.
\item $(2n, \frac{1}{2}+\Z)$. We have to be careful, since $J= 
\begin{pmatrix}
 & 1 \\
 -1 &
\end{pmatrix}$
belongs to $\mr{GSO}_{2n}(\CC)$ if and only if $n$ is even. In that case, the argument proceeds as before, and we find $\omega_v(-1)=-1$. For $n$ odd, we know that $\mr{Cent}_{\mr{GL}_{2n}}(\phi_v(W_{\overline{F}_v}))$ consists of diagonal matrices (this follows from the same calculation as in Lemma \ref{cent}), so if we are to find a $J' \in \mr{GSO}_{2n}(\CC)$ satisfying $J' \phi_v(z) (J')^{-1}= \phi_v(\bar{z})$, it must have the form $J'= Jt$, where $t$ is some diagonal matrix in $\mr{GL}_{2n}(\CC)$. The condition $(J')^2= -1$ is easily seen to imply that $t$ has the form $\diag(t_1, \ldots, t_n, t_1^{-1}, \ldots, t_n^{-1})$, hence lies in $\tT^\vee$. In particular, we cannot choose $J'$ in $\mr{GSO}_{2n}(\CC)$, so in this case we obtain a contradiction: no such $\phi_v \colon W_{F_v} \to \LtG$, and hence no such $\Pi$, can exist.
\item $(2n, \Z)$. Again we may take $J=
\begin{pmatrix}
 & 1_n \\
1_n & 
\end{pmatrix}$, which lies in $\mr{GSO}_{2n}$ exactly when $n$ is even. In this case, $\nu(J)=1$,
and we are forced to have $\omega_v(-1)= 1$. For $n$ odd, an argument as above shows no such parameter, or $\Pi$, can exist.
\end{enumerate}
\endproof
\begin{rmk}
To understand this lemma, and in particular the fact that in the $D_n$ cases we could not have $n$ odd, note that the split real spin group $\mr{Spin}(n, n+1)$ always has discrete series, while the split group $\mr{Spin}(n, n)$ has discrete series if and only if $n$ is even.\footnote{A real orthogonal group $\mr{SO}(p, q)$ is easily seen to have a compact maximal torus if and only if $pq$ is even.} The $\Pi$ we consider will all arise from $\tpi$ such that $\tpi_{\infty}$ is discrete series. Similarly, in the quasi-split non-split case, note that $\mr{Spin}(n-1, n+1)$ has discrete series if and only if $n$ is odd.
\end{rmk}
Finally, we treat the quasi-split but non-split case:
\begin{lemma}
Suppose $\tG$ is the quasi-split group $\mr{GSpin}_{2n}^{\mu}/F$ associated to a non-trivial quadratic extension $F'/F$ cut out by the character $\mu$. Then also in this case $\omega_v(-1)$ is independent of $v \vert \infty$.
\end{lemma}
\proof
For a given $v \vert \infty$, $\tG_v$ may be split or not depending on whether the place(s) of $F'$ above $v$ are real or complex. If real, then $\tG_v$ is split, and the local central character calculations are those of Lemma \ref{splitcase}. If complex, then $\tG_v$ is quasi-split but not split, and we now perform the analogous calculations. First suppose that the $p_{v, i}$ are integers. Writing $\phi_v(j)= (J, c) \in \LtG$ and applying $\iota_{\tG}$ to the relation $\phi_v(jzj^{-1})= \phi_v(\bar{z})$, we obtain
\[
Jh \phi_v(z) h J^{-1}= \phi_v(z)^{-1},
\]
where we abusively write $\phi_v(z)$ both for the element of $\tG^\vee \subset \mr{GL}_{2n}(\CC)$ and for the element of $\LtG$. As before, we find that $Jh$ must have the form $\begin{pmatrix}
 & t \\
 t^{-1} & \end{pmatrix}$ for some $n \times n$ diagonal matrix $t$. Thus $\frac{\det}{\nu^n}(J)= (-1)^{n+1}$, so $J$ can be chosen in $\tG^\vee$ if and only if $n$ is odd, in which case $\nu(J)= 1$.
 
 Similarly, if the $p_{v, i}$ lie in $\frac{1}{2}+\Z$, $Jh$ has the form $\begin{pmatrix}
 & t^{-1} \\
 -t & \end{pmatrix}$, and again $J$ can be taken in $\tG^\vee$ if and only if $\frac{\det}{\nu^n}(J)= (-1)^{n+1}$ equals one, i.e. if and only if $n$ is odd; in this case, $\nu(J)= -1$.
 
To finish the argument, note that $\tG_v$ must be split at all $v \vert \infty$ or non-split at all $v \vert \infty$: if not, then our calculations (based on the existence of a descent $\tpi$ to $\tG$ with regular archimedean L-parameters) show that $n$ must be both odd and even, a contradiction. We have just seen that in the non-split case $\omega_v(-1)$ depends only on whether the $p_{v, i}$ are integers or half-integers, so the proof is complete.
\endproof
\section{Construction of automorphic Galois representations}\label{galois}
We can use the same ideas to `construct' the Galois representations expected to be associated to automorphic representations that are discrete series at infinity on $\mr{GSpin}$ and $\mr{Spin}$ groups, of course building on the deep known results for $\mr{GL}_N$. I stress that the main result, and technical novelty, of this section is Theorem \ref{mixedgalois}, which treats the `mixed-parity' case. Proposition \ref{galoisreps} will be unsurprising, given the calculations of \S \ref{Stwo}. We continue to let $\tG$ denote one of the quasi-split $\mr{GSpin}$ groups $\mr{GSpin}_{2n+1}$, $\mr{GSpin}_{2n}$, $\mr{GSpin}^{\mu}_{2n}$. The calculations of the previous section, in combination with a result of Bella\"{i}che-Chenevier, will allow us to show these Galois representations take values, as hoped, in the appropriate $L$-group. There are a number of more refined statements one could hope for (compare Conjecture 3.2.1, 3.2.2 of \cite{buzzard-gee:alg}), which we can only partially verify: the basic difficulty is that the most precise form of the conjectural relation between automorphic forms and Galois representations is an essentially Tannakian statement, requiring understanding of the Galois representations associated to all functorial transfers of $\tpi$ to general linear groups, whereas we have at our disposal only the single transfer given by $\iota_{\tG}$. We begin with the L-algebraic (on $\tG(\af)$) case, although certainly the more interesting result in this section is for mixed-parity representations (see Theorem \ref{mixedgalois} and Remark \ref{mixedremark} below).
\begin{prop}\label{galoisreps}
Let $F$ be a totally real field, and let $\tpi$ be an L-algebraic,\footnote{By twisting, one can prove a similar result for C-algebraic $\tpi$.} globally generic, cuspidal automorphic representation of $\tG(\af)$ whose archimedean component $\tpi_{\infty}$ is (up to center) a discrete series representation of $\tG(F_{\infty})$. For simplicity, fix an isomorphism $\iota \colon \Qlb \xrightarrow{\sim} \CC$. Then there exists a continuous $\ell$-adic representation
\[
\rho_{\tpi, \iota} \colon \gal{F} \to \LtG(\Qlb),
\]
compatible with the projections to $\gal{F}$, such that $\iota_{\tG} \circ \rho_{\tpi, \iota}$ is the $\ell$-adic representation associated to the functorial transfer $\Pi$ of $\tpi$ to $\mr{GL}_{2n}(\af)$, via the L-homomorphism $\iota_{\tG}$. Moreover, when $\tG= \mr{GSpin}_{2n+1}$, this $\rho_{\tpi, \iota}$ satisfies all parts of Conjecture 3.2.2 of \cite{buzzard-gee:alg}, namely:
\begin{itemize}
\item 
For all finite places $v$ of $F$ outside the finite set $S$ of places where $\tpi$ is ramified, $\rho_{\tpi, \iota}$ is unramified, and $\rho_{\tpi, \iota}(fr_v)$ is $\tG^\vee(\Qlb)$-conjugate to $\iota(\rec_v(\tpi_v)(fr_v))$.
\item For all places $v \vert \ell$ of $F$, $\rho_{\tpi, \iota}|_{\gal{F_v}}$ is de Rham, and it is crystalline if $\tpi_v$ is unramified. For all embeddings $\tau \colon F \into \Qlb$ inducing the place $v$, the $\tau$-labeled Hodge-Tate co-character $\mu_{\tau} \colon \mbb{G}_m \to \tT^\vee$ of $\rho_{\tpi, \iota}|_{\gal{F_v}}$ is (Weyl-conjugate to) the co-character $\mu_{\iota \circ \tau}$ arising from the archimedean L-parameters of $\tpi$: that is, letting $v' \vert \infty$ denote the place induced by $\iota \circ \tau \colon F \into \CC$, $\rec_{v'}(\tpi_{v'})= \phi_{v'}$ takes on $W_{\overline{F}_{v'}}$ the form
\[
\phi_{v'}(z)= (\iota \circ \tau)(z)^{\mu_{\iota \circ \tau}} (\bar{\iota} \circ \tau)(z)^{\mu_{\bar{\iota} \circ \tau}}.
\]
\item For $v \vert \infty$ and $c_v \in \gal{F}$ a complex conjugation at $v$, the $\tG^\vee(\Qlb)$-conjugacy class of $\rho_{\tpi, \iota}(c_v)$ is given by the recipe of \cite[Conjecture 3.2.1]{buzzard-gee:alg}.
\end{itemize}
\end{prop}
\proof
There is an integer $w$ such that $\tpi \otimes |\cdot|^{-w/2}$ is unitary, and it follows that $\omega |\cdot|^{-w}$ is finite-order. Let $\Pi$ denote the transfer of $\tpi$ to $\mr{GL}_{2n}$ provided by Theorem \ref{descent}. By Theorem 4.19 and Corollary 4.24 of \cite{asgari-shahidi:nonsplitGSpin}, $\Pi$ is an isobaric sum $(\sigma_1 \boxplus \cdots \boxplus \sigma_t)  \otimes |\cdot|^{w/2}$ where each $\sigma_i$ is a unitary cuspidal automorphic representation of some $\mr{GL}_{n_i}(\af)$ (where $\sum n_i= 2n$) and satisfies $\sigma_i \cong \sigma_i^{\vee} \otimes \omega |\cdot|^{-w}$. Each $\sigma_i$ is regular, and $\sigma_i |\cdot|^{w/2}$ is L-algebraic, so we associate a Galois representation to $\Pi$, 
\[
\rho_{\Pi, \iota}= \bigoplus_{i=1}^t \rho_{i, \iota} \colon \gal{F} \to \mr{GL}_{2n}(\Qlb),
\]
by applying Theorem 2.1.1 of \cite{blggt:potaut} to each $\sigma_i |\cdot|^{w/2}$. For all $i$, $\rho_{i, \iota}$ preserves a pairing of sign $\omega_{\iota}(c_v)= (-1)^w \omega_v(-1)$, where $c_v$ denotes complex conjugation at $v \vert \infty$, and $\omega_{\iota}$ the associated geometric Galois character: this follows from Theorem \ref{mainthm} above and, more important, Corollary 1.3 of \cite{bellaiche-chenevier:sign}.\footnote{And by the fact that non-self-dual irreducible constituents of $\rho_{\sigma_i, \iota}$ come in (dual) pairs, and on such a pair $r \oplus (r^{\vee} \otimes \omega_{\iota})$ we can put an invariant pairing of any sign we like.} Thus $\rho_{\Pi, \iota}$ preserves (up to scaling) a pairing of sign $(-1)^w \omega_v(-1)$ as well, and the similitude character is just $\omega_{\iota}$. We deduce (from the same calculations as in Lemma \ref{splitcase}) that for $\tG= \mr{GSpin}_{2n+1}$, $\rho_{\Pi, \iota}$ lands in $\mr{GSp}_{2n}(\Qlb)$, while for $\tG= \mr{GSpin}_{2n}$ or $\mr{GSpin}^{\mu}_{2n}$, $\rho_{\Pi, \iota}$ lands in $\mr{GO}_{2n}(\Qlb)$. In the $\mr{GSpin}_{2n}$ case, $\omega_{\Pi}= \omega^{n}$, so $\det \rho_{\Pi, \iota}= \omega_{\iota}^n$, and $\rho_{\Pi, \iota}$ indeed lands in $\mr{GSO}_{2n}(\Qlb)$. Similarly, in the $\mr{GSpin}^{\mu}_{2n}$ case, we claim that $\rho_{\Pi, \iota} \times \mr{id}$ factors through
\[
\iota_{\tG} \colon \mr{GSO}_{2n}(\Qlb) \rtimes \gal{F} \into \mr{GO}_{2n}(\Qlb) \times \gal{F},
\]
but this is immediate from the fact that, for $\gamma \in \gal{F}$, $\gamma|_{F'}$ is trivial if and only if $\mu(\gamma)= \left(\frac{\det}{\nu^n} \circ \rho_{\Pi, \iota}\right)(\gamma)$ is $1$; namely, when trivial we have $\iota_{\tG}(\rho_{\Pi, \iota}(\gamma, \gamma))= (\rho_{\Pi, \iota}(\gamma), \gamma)$, and when non-trivial we have $\iota_{\tG}(\rho_{\Pi, \iota}(\gamma) h, \gamma)= \rho_{\Pi, \iota}(\gamma), \gamma)$. Thus, in all cases we obtain a homomorphism $\rho_{\tpi, \iota} \colon \gal{F} \to \LtG(\Qlb)$.

The fact that, in the $\mr{GSpin}_{2n+1}$ case, we have the more refined $\LtG$-conjugacy of unramified parameters follows immediately from the following elementary observation:
\begin{lemma}[see Lemma 3.6(i) of \cite{larsen:locglobconj}]\label{gspconj}
Two semi-simple elements $x$ and $y$ of $\mr{GSp}_{2n}(\Qlb)$ are $\mr{GSp}_{2n}(\Qlb)$-conjugate if and only if they have the same symplectic multiplier and are $\mr{GL}_{2n}(\Qlb)$-conjugate
\end{lemma}
(Note that this statement fails with $\mr{GSO}_{2n}$ in place of $\mr{GSp}_{2n}$, because of the outer automorphism of $\mr{GSO}_{2n}$ coming from $\mr{GO}_{2n}$; this is why we cannot deduce $\LtG$-conjugacy of unramified parameters in the $\mr{GSpin}_{2n}$ cases.) That $S$ may be taken to be just the set of places at which $\tpi$ is ramified is part of Theorem 4.25 of \cite{asgari-shahidi:nonsplitGSpin}.\footnote{Note that in the case of the quasi-split group $\mr{GSpin}_{2n}^{\mu}$, they lose control at the places ramified in $F'/F$ as well.}

Now let $c_v \in \gal{F}$ be a complex conjugation. When $\tG= \mr{GSpin}_{2n+1}$, we know that $\rho_{\Pi, \iota}(c_v)$ and ${}^t \rho_{\Pi, \iota}(c_v)^{-1} \otimes \omega_{\iota}(c_v)= -{}^t \rho_{\Pi, \iota}(c_v)^{-1}$ are $\mr{GL}_{2n}$-conjugate. Choosing a basis in which $\rho_{\Pi, \iota}(c_v)$ is diagonal, it is then clear that $\rho_{\Pi, \iota}(c_v)$ has eigenvalues $+1$ and $-1$ each with multiplicity $n$. Lemma \ref{gspconj} implies this uniquely determines the $\mr{GSp}_{2n}$-conjugacy class of $\rho_{\tpi, \iota}(c_v) \in \mr{GSp}_{2n}(\Qlb)$; one checks immediately (using the calculation of the Langlands parameter in Lemma \ref{splitcase}) that the conjugacy class predicted by \cite[Conjecture 3.2.1]{buzzard-gee:alg} also satisfies this (determining) property. 

The $\ell$-adic Hodge theory properties follow similarly.
\endproof
\begin{rmk}
It is to be expected that in the cases $\tG= \mr{GSpin}_{2n}$, $\mr{GSpin}^{\mu}_{2n}$, $\rho_{\Pi, \iota}(c_v)$ also has $+1$ and $-1$ eigenspaces each of dimension $n$. \textit{Assuming this}, we can deduce that $\rho_{\tpi, \iota}(c_v)$ is in the $\tG^{\vee}$-conjugacy class predicted by \cite[Conjecture 3.2.1]{buzzard-gee:alg} as long as $\tG_v= \tG \otimes_F F_v$ is split (i.e., in the globally split case, or in the case of quasi-split case when $\mu$ factors through a totally real field). To see this, note that when $\tG_v$ is split, $\mu(c_v)=1$ and $\omega_{\iota}(c_v)=1$, hence $\rho_{\Pi, \iota}(c_v)$ lies in $\mr{SO}_{2n}(\Qlb)$. Since $\rho_{\Pi, \iota}(c_v)$ has eigenvalues equal to $\pm 1$, Lemma 3.7(ii) of \cite{larsen:locglobconj} implies its $\mr{SO}_{2n}$-conjugacy class is uniquely determined by its $\mr{GL}_{2n}$-conjugacy class. It would then follow that $\rho_{\tpi, \iota}(c_v)$ is as in \cite[Conjecture 3.2.1]{buzzard-gee:alg}. If $\tG_v$ is not split, then $\rho_{\Pi, \iota}(c_v)$ lies in $\mr{O}_{2n}-\mr{SO}_{2n}$, but here the corresponding statement that the $\tG^\vee$-conjugacy class of $\rho_{\tpi, \iota}(c_v)$ is determined by the $\mr{GL}_{2n}$-conjugacy class of $\rho_{\Pi, \iota}(c_v)$ fails rather dramatically: up to $\tG^\vee$-conjugacy, $\rho_{\tpi, \iota}(c_v)$ lies in $\tT^\vee \rtimes c$, where $c$ denotes complex conjugation in $\Gal(F'/F)$, and if we write it in the form
\[
\diag(t_1, \ldots, t_n, t_1^{-1}, \ldots, t_n^{-1}),
\]
our knowledge of $\rho_{\Pi, \iota}(c_v)$ tells us that all $t_1, \ldots, t_{n-1}$ are equal to $\pm 1$, with any permutation of these achievable by $\tG^\vee$-conjugacy, but yields no constraint on $t_n$.
\end{rmk}
We can go somewhat farther and associate Galois representations to L-algebraic automorphic representations $\pi$ of $G(\af)$ where $G$ is one of the spin groups $\mr{Spin}_{2n+1}$, $\mr{Spin}_{2n}$, $\mr{Spin}^{\mu}_{2n}$ underlying the corresponding $\tG$. This is not simply a matter of extending $\pi$ to an L-algebraic representation $\tpi$ of $\tG(\af)$, since such an extension need not always exist; see \cite[\S 3.1]{stp:variationsarxiv} for a detailed discussion of such matters. In any case, we would only expect to attach ${}^L G(\Qlb)$-representations (not $\LtG(\Qlb)$) to such a $\pi$. For lack of a suitable generalization of the results of \cite{bellaiche-chenevier:sign}, we will only produce the $\mr{PGL}_{2n}(\Qlb)$-representations corresponding to the `projectivization' of $\iota_{\tG}$, which gives an embedding $\iota_G \colon {}^L G \into \mr{PGL}_{2n} \times \gal{F}$. We will have to make use of the more general construction of automorphic Galois representations over CM fields (again, see \cite[Theorem 2.1.1]{blggt:potaut} for a precise statement); for a CM field $L$ with totally real subfield $F$, we call an automorphic representation $\Pi$ of $\mr{GL}_N(\mbf{A}_L)$ polarizable if $\Pi^c \cong \Pi^{\vee} \otimes \mr{BC}_{L/F}(\omega)$, where $\omega \colon C_F \to \CC^\times$ is a Hecke character of $F$ such that $\omega_v(-1)$ is independent of $v \vert \infty$. 
\begin{thm}\label{mixedgalois}
Let $F$ be a totally real field, and let $\pi$ be a globally generic, L-algebraic discrete series at infinity, cuspidal automorphic representation of $G(\af)$. Then there exist continuous $\ell$-adic representations
\[
\rho_{\pi, \iota} \colon \gal{F} \to \mr{PGL}_{2n}(\Qlb)
\]
such that
\begin{enumerate}
\item For all finite places $v$ of $F$ outside the finite set $S$ of places where $\pi$ is ramified, $\rho_{\pi, \iota}$ is unramified, and $\rho_{\pi, \iota}(fr_v)$ is $\mr{PGL}_{2n}(\Qlb)$-conjugate to $\iota(\iota_G \circ \rec_v(\pi_v)(fr_v))$.
\item For all $v \vert \ell$, $\rho_{\pi, \iota}$ is de Rham, and it is crystalline if $\pi_v$ is unramified. Its Hodge-Tate co-characters  are given by the same archimedean recipe, but then composed with $\iota_G$, as in Proposition \ref{galoisreps}.
\end{enumerate}
\end{thm}
\proof
As in \cite[Proposition 3.1.14]{stp:variationsarxiv}, we may choose a W-algebraic extension $\tpi$ of $\pi$ to $\tG(\af)$,\footnote{Unfortunately, in \cite[\S 3.1]{stp:variationsarxiv} the extension of automorphic representations from $G(\af)$ to $\tG(\af)$ was written assuming the center $Z_{\tG}$ of $\tG$ was a torus; this does not hold for $\tG= \mr{GSpin}_{2n}$, but all that is in fact required is that the quotient $Z_{\tG}/Z_G$ be a torus.} and then apply Theorem \ref{descent} to obtain the transfer $\Pi$ of $\tpi$ to $\mr{GL}_{2n}(\af)$. Twisting by type A Hecke characters, we will be able to associate $\mr{GL}_{2n}$-valued Galois representations over (quadratic) CM extensions $L/F$. The main obstacle to descent to a projective $\gal{F}$-representation is reducibility of these $\gal{L}$-representations; this is surmounted by the patching lemma, originally due to Blasius-Rogawski (\cite[Proposition 4.3]{blasius-rogawski:motiveshmfs}; see too the more general and very clear presentation in \cite{sorensen:patching}). Now, the extension $\tpi$ of \cite[Proposition 3.1.14]{stp:variationsarxiv} has, by construction, finite-order central character. From this and the description of discrete series L-parameters, we see that $\Pi_v$ is L-algebraic for certain $v \vert \infty$ and C-algebraic for others, yielding a partition $S_{\infty}= S_L \sqcup S_C$ of the archimedean places of $F$.\footnote{To be explicit in a particular case, suppose $G= \mr{Spin}_{2n+1}$. $G$ is simply-connected, so $\rho$ is in the weight lattice, and our discrete series L-parameter at $v \vert \infty$ is determined by a single element $\mu_v \in \rho+X^\bullet(T)= X^\bullet(T)$, for $T$ a maximal torus. Letting $\tT \supset T$ be the corresponding maximal torus of $\tG$, we have, in suitable coordinates, a Cartesian diagram 
\[
\xymatrix{
X^\bullet(\tT)= \oplus_{i=1}^n \Z \chi_i \oplus \Z(\chi_0+ \frac{\sum_1^n \chi_i}{2}) \ar[r] \ar[d] & X^\bullet(Z_{\tG}) \ar[d] \\
X^\bullet(T)= \oplus_{i=1}^n \Z \chi_i \oplus \Z(\frac{\sum_1^n \chi_i}{2}) \ar[r] & X^\bullet(Z_G).
}
\]
The lifted L-parameter on $\tG(F_v)$ is given by $\tilde{\mu}_v \in \frac{1}{2} X^\bullet(\tT)$ projecting to zero (since the central character of $\tpi$ is finite-order) in $X^\bullet(Z_{\tG})$ and to $\mu_v \in X^\bullet(T)$. In particular, if $\mu_v$ lies in $\oplus \Z \chi_i$, then $\Pi_v$ is L-algebraic; and if $\mu_v$ lies in $\frac{\sum \chi_i}{2} + \oplus \Z \chi_i$, then $\Pi_v$ is C-algebraic.
} 
For simplicity, for each $v \vert \infty$ of $F$, we fix an embedding $\sigma_v \colon \overline{F} \into \CC$ whose restriction to $F$ induces $v$; this will help us compare infinity-types of Hecke characters of varying CM extensions of $F$. We consider almost all quadratic CM extensions $L/F$ of the form $L= F(\sqrt{-p})$, where $p$ is a rational prime; we can throw away any finite number, and it will be easiest to think about some of the arguments below by throwing out all $F(\sqrt{-p})$ such that $F$ is ramified at some place above $p$. Let us denote the set of such $L$ by $\mc{I}$: it has the property (called `strongly $\emptyset$-general in \cite{sorensen:patching}) that for any finite set $\Sigma$ of places of $F$, there is an $L \in \mc{I}$ in which every $v \in \Sigma$ splits completely. For each $v \in S_{\infty}$, our fixed $\sigma_v$ induces $\sigma_w \colon L_w \xrightarrow{\sim} \CC$ for the unique infinite place $w$ of $L$ above $v$, and with reference to these embeddings we select an infinity-type for a Hecke character $\psi_L$ of $L$, letting 
\[
\psi_{L_w}(z)=
\begin{cases}
1 & \text{if $w|_F \in S_L$,}\\
\sigma_w(z)/|\sigma_w(z)| & \text{if $w|_F \in S_C$.}
\end{cases}
\]
By \cite{weil:characters}, there indeed exists a Hecke character $\psi_L \colon C_L \to \CC^\times$ with this infinity-type. For all $L \in \mc{I}$, we can then form the L-algebraic automorphic representation $\mr{BC}_{L/F}(\Pi) \otimes \psi_L$. It need not be cuspidal, but recall that $\Pi$ is isobaric with all cuspidal constituents $\sigma_i$ satisfying $\sigma \cong \sigma^\vee \otimes \tomega$, where $\tomega$ is the central character of $\tpi$. Throwing out the finite number of $L \in \mc{I}$ such that $\mr{BC}_{L/F}(\sigma_i)$ is non-cuspidal for some $i$ (the resulting set remains `strongly $\emptyset$-general'), each $\mr{BC}_{L/F}(\sigma_i) \otimes \psi_L$ is a polarizable, regular, algebraic, cuspidal automorphic representation of $\mr{GL}_{n_i}(\mbf{A}_L)$; this follows since
\[
(\mr{BC}_{L/F}(\sigma_i) \otimes \psi_L)^{c} \cong \mr{BC}_{L/F}(\sigma_i) \otimes \psi_L^c \cong (\mr{BC}_{L/F}(\sigma_i) \otimes \psi_L)^\vee \otimes (\tomega \psi_L^{1+c}),
\]
where $c$ denotes a complex conjugation. For polarizability, note that $\psi_L^{1+c}$ descends to a Hecke character of $F$, and that, while $\tomega_v(-1)$ is \textit{not} independent of $v \vert \infty$, the infinity-type of $\psi_L$ precisely ensures (by the central character calculations of \S \ref{Stwo}) that $(\tomega \psi_L^{1+c})_v(-1)$ is independent of $v \vert \infty$ (for either descent of this Hecke character to $F$). By Theorem 2.1.1 of \cite{blggt:potaut}, we can therefore associate a compatible system of $\ell$-adic representations 
\[
\rho_{\mr{BC}_{L/F}(\Pi) \otimes \psi_L, \ell} \colon \gal{L} \to \mr{GL}_{2n}(\Qlb)
\]
compatible with the local parameters of $\mr{BC}_{L/F}(\Pi) \otimes \psi_L$. Although $\psi_L$ is not type $A_0$, and consequently does not have an associated $\ell$-adic Galois representation, we can find a Galois character $\hat{\psi}_L \colon \gal{L} \to \Qlb^\times$ with the property that $\hat{\psi}_L^{c-1}$ is the (geometric) Galois character associated to $\psi_L^{c-1}$, and, moreover, $\psi_L^2$ is type $A_0$ with associated geometric Galois character differing from $(\hat{\psi}_L)^2$ by a finite-order character of $\gal{L}$.\footnote{To see this, first let $\hat{\psi} \colon \gal{F} \to \Qlb^\times$ be a Galois character such that $(\psi^2)_{\iota}\cdot \hat{\psi}^{-2}$ is finite-order-- here we write $(\cdot)_{\iota}$ for the Galois character associated to a type $A_0$ Hecke character via $\iota$-- and therefore $(\psi^{c-1})_{\iota}$ and $\hat{\psi}^{c-1}$ also differ by a finite-order character. Invoking \cite[Lemma 3.3.4]{stp:variationsarxiv}, we can find a type $A$ Hecke character $\psi_1$ of $L$ such that $(\psi_1^{c-1})_{\iota}= \hat{\psi}^{c-1}$, and it is easy to see (by checking the infinity-type) that $\psi_1/\psi$ is finite-order, hence has an associated Galois character $\gamma \colon \gal{F} \to \Qlb$. Clearly $(\hat{\psi} \gamma^{-1})^{c-1}= (\psi^{c-1})_{\iota}$, so the claim is proven.} The $\ell$-adic representations
\[
\rho_{L, \ell}:=  \rho_{\mr{BC}_{L/F}(\Pi) \otimes \psi_L, \ell} \otimes \hat{\psi}_L^{-1} \colon \gal{L} \to \mr{GL}_{2n}(\Qlb)
\]
are then $\Gal(L/F)$-invariant.

We now fix an $\ell$ and a base-point $L_0 \in \mc{I}$, and we compare the various $\rho_L= \rho_{L, \ell}$: for $L, L' \in \mc{I}$, we find that as $\gal{LL'}$-representations
\[
\rho_{L} \cong \rho_{L'} \otimes \Psi(L, L'),
\] 
where $\Psi(L, L')$ is the (finite-order) Galois character of $\gal{LL'}$ given by
\[
\Psi(L, L')= \left( \frac{\psi_L}{\psi_{L'}} \right)_{\ell} \cdot \frac{\hat{\psi}_{L'}}{\hat{\psi}_L},
\]
where we abusively let $\left( \frac{\psi_L}{\psi_{L'}} \right)_{\ell}$ denote the $\ell$-adic character associated to the finite-order (by our choices of infinity-types) Hecke character $\frac{\psi_L}{\psi_{L'}} \colon C_{LL'} \to \CC^\times$. We clearly have the `co-cycle relation' 
\[
\Psi(L, L') \cdot \Psi(L', L'')= \Psi(L, L'')
\]
on the triple intersections $\gal{LL'L''}$. Each character $\Psi(L, L')$ is $\Gal(LL'/F)$-invariant (using the relations $\psi_L^{c-1}= \hat{\psi}_L^{c-1}$), and we would like to descend this co-cycle relation, but invariant Hecke characters of a number field need not descend through non-cyclic (even biquadratic, as here) extensions.\footnote{For example, let $F \subset K \subset L$ be a tower of number fields with $\Gal(L/F)$ the dihedral group with $8$ elements, and $\Gal(L/K)$ its central $\Z/2$; taking $\psi$ to be the Hecke character of $K$ cutting out the extension $L/K$, $\psi$ is $\Gal(K/F)$-invariant, but it clearly does not descend.} Nevertheless, we can get away with somewhat less. Fix an $L_1 \in \mc{I}$ linearly disjoint from $L_0$ over $F$, and fix a descent $\Psi(L_0, L_1)_{L_1}$ of $\Psi(L_0, L_1)$ to a character of $\gal{L_1}$; this descent is determined up to $\Gal(L_0/F) \cong \Gal(L_0L_1/L_1)$-twist. For any other $L \in \mc{I}$, linearly disjoint from $L_0L_1$ over $F$, we define $\Psi(L_0, L)_L$ to be the unique descent of $\Psi(L_0, L)$ to a character of $\gal{L}$ with the property that
\[
\left( \Psi(L_0, L_1)_{L_1} \Psi(L_1, L) \Psi(L_0, L)_L^{-1} \right)|_{\gal{L_1L}} = 1.
\]
(This triple-product is trivial on $\gal{L_0L_1L}$ by the co-cycle relation, and so it is trivial on $\gal{L_1L}$ for exactly one of the descents of $\Psi(L_0, L)$ to $\gal{L}$: note that $L_0L \neq L_1L$, so $L_0L_1L/L$ is in fact biquadratic.) 

Next, for all $L \neq L_0$, we replace $\rho_L$ by its twist
\[
\rho^{\star}_L= \rho_L \otimes \Psi(L_0, L)_L \colon \gal{L} \to \mr{GL}_{2n}(\Qlb).
\]
We claim that for any two $L, L' \in \mc{I}$, $\rho_L^{\star}|_{\gal{LL'}} \cong \rho_{L'}^{\star}|_{\gal{LL'}}$. We have
\begin{align*}
\rho_L^{\star}|_{\gal{LL'}}&= (\Psi(L_0, L)_L \otimes \rho_L)|_{\gal{LL'}}= (\Psi(L_0, L)_L \Psi(L, L') \otimes \rho_{L'} )|_{\gal{LL'}} \\
&= (\Psi(L_0, L)_L \Psi(L, L') \Psi(L_0, L')_{L'}^{-1})_{\gal{LL'}} \otimes \rho_{L'}^{\star}|_{\gal{LL'}},
\end{align*}
where note that the twisting character is trivial after further restriction to $\gal{L_0LL'}$. To see that it is in fact trivial on $\gal{LL'}$, consider the three expressions
\begin{align*}
A&= \left( \Psi(L_0, L_1)_{L_1} \Psi(L_1, L) \Psi(L_0, L)_L^{-1} \right)|_{\gal{L_1L}} = 1; \\
B&= \left( \Psi(L_0, L_1)_{L_1} \Psi(L_1, L') \Psi(L_0, L')_{L'}^{-1} \right)|_{\gal{L_1L'}} = 1; \\
C&= \left( \Psi(L_0, L)_L \Psi(L, L') \Psi(L_0, L')_{L'}^{-1} \right)_{\gal{LL'}} = \text{?}
\end{align*}
On the triple intersection, we have (by canceling and applying the co-cycle relation) $AB^{-1}C|_{\gal{L_1LL'}}=1$, hence $C|_{\gal{L_1LL'}}=1$. But we also know that $C|_{\gal{L_0LL'}} =1$, so, since $L_0LL' \neq L_1LL'$ (by our choice of $\mc{I}$, these fields are not ramified at the same set of primes), $C$ itself equals $1$, and we conclude that $\rho_{L}^{\star}|_{\gal{LL'}} \cong \rho_{L'}^{\star}|_{\gal{LL'}}$.

Finally, we can apply the patching lemma (\cite[Lemma 1]{sorensen:patching}) to deduce that there is a continuous semi-simple representation $\rho \colon \gal{F} \to \mr{GL}_{2n}(\Qlb)$ such that $\rho|_{\gal{L}} \cong \rho^{\star}_L$ for all $L \in \mc{I}$. The projectivization of this $\rho$ is the desired $\rho_{\pi, \iota}$ as in the statement of our proposition: for all $L$, and for all places $v$ of $F$ split in $L$ (say $v= w \bar{w}$), $\rho(fr_v)$ has the same projective parameter as $\rho_{\BC_{L/F}(\Pi) \otimes \psi_L}(fr_w)$, hence the same projective parameter as $\Pi_v$. Varying $L$, since $\mc{I}$ is strongly $\emptyset$-general, we obtain the desired compatibility at almost all places $v$ of $F$.
\endproof
\begin{rmk}\label{mixedremark}
\begin{itemize}
\item The same argument clearly yields associated projective $\ell$-adic representations for any mixed-parity, regular, essentially self-dual cuspidal automorphic representation of $\mr{GL}_{2n}(\af)$; see Corollary \ref{pm1liftstotreal}.
\item Note that, in contrast to Proposition \ref{galoisreps}, we lose control of the explicit set of places of $F$ where we do not know local-global compatibility: this occurs because in extending $\pi$ to the desired $\tpi$ on $\tG(\af)$, we might have to allow ramification at additional primes, at which we cannot apply Theorem \ref{descent}. This does not interfere with the claim that $\rho_{\pi, \iota}|_{\gal{F_v}}$ is crystalline at $v \vert \ell$ when $\pi_v$ is unramified, because we can choose the extension $\tpi$ to be unramified at any finite set of primes at which $\pi$ is unramified.
\end{itemize}
\end{rmk}

\section{A sharper result for mixed-parity automorphic representations}\label{sectionmixed}
Claus Sorensen has suggested to me that in situations such as Theorem \ref{mixedgalois}, in which one works with W-algebraic, but \textit{mixed-parity} (in the terminology of Definition \ref{mixedparitydef}), automorphic representations $\Pi$ of $\mr{GL}_N(\af)$,\footnote{Note that this rules out certain W-algebraic representations that are strange hybrids of L- and C-algebraic; see \cite[Example 2.5.6]{stp:variationsarxiv} for some discussion of this.} one should be able to find something stronger than an associated projective Galois representation, but rather a $\mr{GL}_{N}(\Qlb)/\{\pm 1\}$-valued representation. In this section, we provide some illustrations of Sorensen's idea: in the first part, we review and refine certain lifting results from \cite{stp:variationsarxiv}, and in what follows we give a detailed example. First, however, we make a further remark on the mixed-parity condition. As we have seen, for any even $N$ there exist such mixed-parity representations over totally real fields; the situation is markedly different for odd $N$.\footnote{Over totally real fields; of course for any $N$ there are mixed-parity representations over CM fields, simply because there are type $A$, but not $A_0$, Hecke characters} It is well-known that no mixed-parity Hecke characters ($N=1$) exist for totally real fields, and this in fact continues to be the case for any odd $N$:
\begin{lemma}\label{nomixedodd}
Let $F$ be a totally real number field, and let $\Pi$ be a cuspidal automorphic representation of $\mr{GL}_N(\af)$, for an odd positive integer $N$. Then $\Pi$ cannot be mixed-parity.
\end{lemma}
\proof
The key point is that, although $\Pi$ need not be essentially self-dual, its infinity-type is; this follows from an analysis of the possible representations of the Weil group $W_{\RR}$. We may assume $\Pi$ is unitary, and then from Clozel's archimedean purity lemma (\cite[Lemme 4.9]{clozel:alg}) and the structure of $W_{\RR}$ it follows that $\rec_v(\Pi_v)|_{W_{\overline{F}_v}}$ takes the form
\[
z \mapsto \oplus_{i=1}^N \left( \frac{z}{|z|} \right)^{m_{i, v}},
\]
where, ordering the $m_{i, v}$ so that $m_{1, v} \geq m_{2, v} \geq \dots \geq m_{N, v}$, we have
\[
m_{i, v}= -m_{N+1-i, v}
\]
for all $i$. In particular, $m_{\frac{N+1}{2}, v}=0$, and hence all $m_{i, v}$ are even (by our definition of mixed-parity, the $m_{i, v}$ are either all odd or all even for a given $v$). This holds for all $v$, so $\Pi$ (with this unitary normalization) is L-algebraic.
\endproof
\begin{rmk}
We can also regard this lemma as an automorphic instance of `Hodge symmetry.' Such a mixed-parity $\Pi$ would be L-algebraic on $\mr{SL}_N(\af)$, hence would be expected to yield geometric, even motivic, $\mr{PGL}_N(\Qlb)$-valued representations of $\gal{F}$. Precisely when $N$ is odd, Hodge symmetry for motivic Galois representations (see \cite[\S 3.2]{stp:variationsarxiv}) forces the existence of \textit{geometric} lifts to $\mr{GL}_N(\Qlb)$ of these projective representations. Such a geometric lift in turn should correspond to an L-algebraic automorphic representation $\widetilde{\Pi}$ that is-- almost!-- twist-equivalent to $\Pi$, although here there is a subtlety arising frome endoscopic phenomena: strictly speaking, we deduce from compatibility of the local parameters only the existence of a character $\chi \colon \af^\times \to \CC^\times$, not necessarily factoring through a Hecke character, such that $\widetilde{\Pi} \cong \Pi \otimes \chi$. But some power of $\chi$ will be a type $A_0$ Hecke character  (compare \cite[Lemma 3.1.9]{stp:variationsarxiv}), so ($F$ is totally real) the infinity-type of $\chi$ is still constrained to that of a rational power of the absolute value. The isomorphism $\widetilde{\Pi} \cong \Pi \otimes \chi$ is then easily seen to contradict the assumption that $\Pi$ was mixed-parity.
\end{rmk}

\subsection{General Galois lifting results}
We take the occasion to streamline, and slightly generalize, some of the arguments of \cite[\S 3.2]{stp:variationsarxiv}. Recall the basic problem, as posed by Brian Conrad in \cite{conrad:dualGW}: given a quotient $H' \onto H$ of linear algebraic groups over $\Qlb$ with kernel equal to a central torus in $H'$, and given a geometric (in the sense of Fontaine-Mazur) Galois representation $\rho \colon \gal{F} \to H(\Qlb)$, when does there exist a geometric lift $\tilde{\rho}$ completing the diagram
\[
\xymatrix{
& H'(\Qlb) \ar[d] \\
\gal{F} \ar[r]_{\rho} \ar@{-->}[ur]^{\tilde{\rho}} & H(\Qlb)
}
\]
It turns out that we can answer this question if $\rho$, in addition to being de Rham at places above $\ell$, moreover satisfies a `Hodge-symmetry' hypothesis that will always hold under the Fontaine-Mazur and Tate conjectures.\footnote{Essentially the same arguments should yield `if and only if' statements that don't demand any faith in these deep conjectures to be convincing; see Proposition \ref{reallift} below, where we obtain such an `if and only if' statement for totally real $F$.} Here we focus on the case in which $F$ admits a real embedding; this was not treated in full generality in \cite{stp:variationsarxiv}, although the argument is in fact easier than the totally imaginary case; what is harder is pinning down the correct hypotheses.

The analysis of the lifting problem is ultimately reduced to the case of connected reductive $H$ and $H'$. For this reason, and to make the automorphic analogy plain, we will for the time being replace $H'$ and $H$ by the dual groups $\tG^\vee$ and $G^\vee$ of connected reductive groups $G$ and $\tG$ (which in the automorphic context would live over $F$, but for us this is irrelevant; one can think of them as split groups over $F$, or just as groups over some algebraically closed field of characteristic zero) constructed as follows: $G$ is arbitrary, and $\tG$ is formed by extending the center $Z_G$ of $G$ to a multiplicative type group $\tZ$ with the property that $S= \tZ/Z_G$ is a \textit{torus}; that is, we set
\[
\tG= (G \times \tZ)/{Z_G},
\]
where $Z_G$ is embedded anti-diagonally. Equivalently, we can fix $\tG$ and let $G$ be any connected subgroup containing the derived group of $\tG$; the quotient is a torus since $\tG$ is connected. The fact that $S$ is a torus is equivalent to the kernel of $\tG^\vee \onto G^\vee$ being a (central) torus.\footnote{Note that $\tZ$ itself need not be a torus; for instance, we could have $G= \mr{Spin}_{2n}$, $\tG= \mr{GSpin}_{2n}$.} Let $T$ be a maximal torus in $G$, with dual maximal torus $T^\vee$ in $G^\vee$; then $\tT= (T \times \tZ)/Z_G$ is a maximal torus in $\tG$, with dual $\tT^\vee$ in $\tG^\vee$. The relevant group theory will be encoded by the following two (exact) diagrams:
\[
\xymatrix{
0 \ar[r] & X^\bullet(T^\vee) \ar[r] & X^\bullet(\tT^\vee) \ar[r] & X^\bullet(S^\vee) \ar[r] & 0,
}
\]
and, writing, for an abelian group $X$, $X_{tor}$ for its torsion subgroup, and $X_{tor}^D= \Hom(X_{tor}, \Q/\Z)$ for the Pontryagin dual of $X_{tor}$,
\[
\xymatrix{
X^\bullet(T^\vee) \oplus X^\bullet((\tZ^0)^\vee) \ar[r] & X^\bullet(\tT^\vee) \ar@{-->}[d] \ar[r] & \Ext^1(X^\bullet(Z_G), \Z) \ar[r] & \Ext^1(X^\bullet(\tZ), \Z) \ar[r] & 0 \\
0 \ar[r] & K \ar[r] & X^\bullet(Z_G)_{tor}^D \ar[u]^{\sim} \ar[r] & X^\bullet(\tZ)_{tor}^D \ar[u]^{\sim} \ar[r] & 0;
}
\]
here $K$ is by definition the kernel of the bottom right map.
The top sequence comes from taking character groups and then applying $\Hom( \cdot, \Z)$ to the sequence
\[
1 \to Z_G \to T \times \tZ \to \tT \to 1,
\]
noting that 
\[
X^\bullet((\tZ^0)^\vee) \cong X_{\bullet}(\tZ^0) \xrightarrow{\sim} \Hom(X^\bullet(\tZ), \Z).
\] 
Note too that 
\[
X^\bullet((\tZ^0)^\vee) \xrightarrow{\sim} X^\bullet(\tG^\vee).
\] 
The vertical isomorphisms are boundary maps associated to the sequence 
\[
0 \to \Z \to \Q \to \Q/\Z \to 0,
\]
and $K$ maps isomorphically (via the arrows in the diagram) onto the image of $X^\bullet(\tT^\vee)$ in $\Ext^1(X^\bullet(Z_G), \Z)$. Note that the Pontryagin dual of $K$ admits a canonical isomorphism 
\begin{equation}\label{grouptheory}
\coker \left( X^\bullet(\tZ)_{tor} \to X^\bullet(Z_G)_{tor} \right) \xrightarrow{\sim} K^D.
\end{equation}
With this group theory in mind, our lifting result boils down to a lemma on Galois characters and (only necessary in the totally imaginary case) a simple trick in the representation theory of $\tG^\vee$. For the lemma, and subsequently, we will use the language of `Hodge-Tate quasi-cocharacters' for Galois representations that are not necessarily Hodge-Tate but do have rational Hodge-Tate-Sen weights; there are various ways of making this precise-- in \cite{stp:variationsarxiv} we used the much more general Sen theory, although here the usual Hodge-Tate theory plus a little group theory would suffice.
\begin{lemma}\label{galoischaracters}
Let $S^\vee$ be a $\Qlb$-torus, as above. Suppose we are given a collection, indexed by embeddings $\tau \colon F \into \Qlb$, $\{\mu_{\tau}\}_{\tau}$ of quasi-cocharacters of $S^\vee$, i.e. elements of $X_{\bullet}(S^\vee)_{\Q}= X_{\bullet}(S^\vee) \otimes_{\Z} \Q$. Then there exists a Galois character $\psi \colon \gal{F} \to S^\vee(\Qlb)$ with $\tau$-labeled Hodge-Tate quasi-cocharacters $\mu_{\tau}$ if and only if:
\begin{enumerate}
\item if $F$ has a real embedding, the $\mu_{\tau} \in X_{\bullet}(S^\vee)_{\Q}$ are independent of $\tau$;
\item if $F$ is totally imaginary with maximal CM subfield $F_{cm}$, the $\mu_{\tau}$ depend only on $\tau_0= \tau|_{F_{cm}}$ and $\mu_{\tau_0}+ \mu_{\tau_0 \circ c}$ is independent of $\tau_0$.
\end{enumerate}
\end{lemma} 
\proof
Choosing coordinates, i.e. an isomorphism $S^\vee \cong \mbf{G}_m^r$ for some $r$, this becomes an existence problem for $r$ different Galois characters $\gal{F} \to \Qlb^\times$; the result follows from \cite[Corollary 2.3.16]{stp:variationsarxiv} (the imaginary case) and \cite[Lemma 2.3.17]{stp:variationsarxiv} (the real case; note that if $F$ has a real embedding, its maximal subfield with well-defined complex conjugation must be totally real, so the type A Hecke characters of $F$ are just finite-order twists of rational powers of the absolute value).
\endproof
Let us return to our geometric $\rho \colon \gal{F} \to G^\vee(\Qlb)$. Let $\mu_{\tau} \in X_{\bullet}(T^\vee)$ denote its $\tau$-labeled Hodge-Tate cocharacters (defined up to Weyl-conjugacy, although we fix representatives). In proving the existence of some lift $\tilde{\rho}$ of $\rho$ to $\tG^\vee$, Conrad (\cite[Proposition 5.3]{conrad:dualGW})) exploits the existence of an isogeny-complement $\tH$ to $S^\vee$ in $\tG^\vee$; that is, a closed subgroup $\tH \subset \tG^\vee$ such that $\tH \to G^\vee$ is surjective with finite kernel. In our context, we may assume $\tH$ is connected reductive, and Tate's theorem ($H^2(\gal{F}, \Q/\Z)=0$) ensures that $\rho$ lifts to a homomorphism valued in a subgroup of $\tG^\vee$ whose neutral component is $\tH$. 

Let us write $T_{\tH}$ for the maximal torus of $\tH$ lying above $T^\vee$. The homomorphism $X_{\bullet}(T_{\tH}) \to X_{\bullet}(T^\vee)$ is not surjective, but it induces an isomorphism of quasi-cocharacter groups:
\[
X_{\bullet}(T_{\tH})_{\Q} \xrightarrow{\sim} X_{\bullet}(T^\vee)_{\Q}.
\]
We can therefore uniquely lift each $\mu_{\tau}$ to an element $\tmu_{\tau}$ of $X_{\bullet}(T_{\tH})_{\Q}$.
\begin{lemma}
There exists a lift $\tilde{\rho} \colon \gal{F} \to \tG^\vee(\Qlb)$ of $\rho$ with labeled Hodge-Tate quasi-cocharacters given by the images of $\tmu_{\tau}$ (which we still denote $\tmu_{\tau}$) under $X_{\bullet}(T_{\tH})_{\Q} \to X_{\bullet}(\tT^\vee)_{\Q}$.
\end{lemma}
\proof
Indeed, the lift described in the preceding paragraphs, landing in some finite-index supergroup of $\tH$, has the desired property.
\endproof
These lifts define pairings $\langle \tmu_{\tau}, \cdot \rangle \colon X^\bullet(\tT^\vee) \to \Q$, which actually factor through pairings (compare \cite[Lemma 3.2.2]{stp:variationsarxiv})
\[
\langle \tmu_{\tau}, \cdot \rangle \colon \frac{X^\bullet(\tT^\vee)}{X^\bullet(T^\vee) \oplus X^\bullet(\tG^\vee)} \to \Q/\Z,
\]
and can therefore be identified with elements 
\[
\theta_{\rho, \tau} \in \coker \left( X^\bullet(\tZ)_{tor} \to X^\bullet(Z_G)_{tor} \right).
\]
(Recall equation \ref{grouptheory}.) Even if $\rho$ takes values in a non-connected group $H$, we can still define these $\theta_{\rho, \tau}$ since the Hodge-Tate co-characters are valued in the identity component $H^0$. Now we describe the lifting obstruction when $F$ admits a real embedding:
\begin{prop}\label{reallift}
Let $H' \onto H$ be a surjection of linear algebraic groups with central torus kernel, and write the identity components $(H')^0 = \tG^\vee$ and $H^0 = G^\vee$ where $G$ and $\tG$ are connected reductive groups as above. Suppose $F$ is a number field with at least one real embedding, and that $\rho \colon \gal{F} \to H(\Qlb)$ is a geometric Galois representation. Then $\rho$ admits a lift to a geometric representation $\tilde{\rho} \colon \gal{F} \to H'(\Qlb)$ if, and only if, the elements $\theta_{\rho, \tau} \in \coker \left( X^\bullet(\tZ)_{tor} \to X^\bullet(Z_G)_{tor} \right)$ are independent of $\tau \colon F \into \Qlb$.
\end{prop}
\proof
We first assume $\rho$ has connected algebraic monodromy group, in which case we can replace $H$ by $G^\vee$ and $H'$ by $\tG^\vee$. The elements $\theta_{\rho, \tau}$ are independent of $\tau$ if and only if the pairings 
\[
\langle \tmu_{\tau}, \cdot \rangle \colon X^\bullet(S^\vee) \cong X^\bullet(\tT^\vee)/X^\bullet(T^\vee) \to \Q/\Z
\]
are independent of $\tau$. The lift $\tilde{\rho}$ chosen above can be twisted (via a character $\gal{F} \to S^\vee(\Qlb)$) to a geometric lift of $\rho$ if and only this induced collection of pairings $X^\bullet(S^\vee) \to \Q/\Z$ arises from a collection $\lambda_{\tau} \in X_{\bullet}(S^\vee)_{\Q}$\footnote{Under the surjection $X_{\bullet}(S^\vee)_{\Q} \onto X_{\bullet}(S^\vee) \otimes_{\Z} \Q/\Z$} equal to the $\tau$-labeled quasi-cocharacters of some Galois character $\psi \colon \gal{F} \to S^\vee(\Qlb)$. The proposition follows from Lemma \ref{galoischaracters}.

For general (not necessarily connected) $\rho$, we still choose the initial lift $\tilde{\rho}$ as above. Even though $\tilde{\rho}$ does not land in $\tG^\vee$, that it twists to a geometric lift is still equivalent to the existence of a Galois character $\psi \colon \gal{F} \to S^\vee(\Qlb)$ with the prescribed Hodge-Tate quasi-cocharacters. That the $\theta_{\rho, \tau}$ are independent of $\tau \colon F \into \Qlb$ is equivalent to the existence of such a Galois character over $F$. 
\endproof

We will spend some time unpacking this statement. First, for completeness, and contrast, we recall without proof the situation in the totally imaginary case:
\begin{thm}[Theorem 3.2.10 of \cite{stp:variationsarxiv}]
Let $H' \onto H$ be an arbitrary surjection of linear algebraic groups over $\Qlb$ with kernel equal to a central torus in $H'$. Let $F$ be totally imaginary, and let $\rho \colon \gal{F} \to H(\Qlb)$ be a geometric Galois representation. Then there exists a geometric lift $\tilde{\rho} \colon \gal{F} \to H'(\Qlb)$ of $\rho$ as long as $\rho$ satisfies `Hodge-symmetry' (\cite[Hypothesis 3.2.4]{stp:variationsarxiv}).
\end{thm}

Now let us return to the totally real case with some examples. The case in which $\tZ$ is a torus, in which case the elements $\theta_{\rho, \tau}$ are simply elements of $X^\bullet(Z_G)_{tor}$, was treated, with examples, in \cite{stp:variationsarxiv}; there (Corollary 3.2.8) we also explained, under the `Hodge-symmetry' hypothesis, how far these elements could be from being independent of $\tau$ (namely, that it reduces to a parity obstruction: only the two-torsion in $X^\bullet(Z_G)$ matters). Sorensen has suggested some interesting examples of $\mr{GL}_N/\mu_m \to \mr{PGL}_N$ lifting problems; we first digress to explain what this looks like on the dual (`automorphic') side.
\subsection{The groups $m\mr{GL}_N$}
The group $\mr{GL}_N$ has fundamental group $\Z$, so there is a unique $\Z/m$-covering space for all natural numbers $m$.
\begin{lemma}
Let $m$ and $N$ be natural numbers. Denote by $m\mr{GL}_N$ the unique connected $\Z/m$-cover of $\mr{GL}_N$, i.e. the Cartesian product
\[
\xymatrix{
m\mr{GL}_N \ar[r] \ar[d] & \mbf{G}_m \ar[d]^{[m]} \\
\mr{GL}_N \ar[r]_{\det} & \mbf{G}_m.
}
\]
Then the dual group $(m\mr{GL}_N)^\vee$ is isomorphic to $\mr{GL}_N/\mu_m$, the quotient by the (central) $m^{th}$ roots of unity.  
\end{lemma} 
\proof
Omitted. To set up the calculation, take a maximal torus $T_N$ of $\mr{GL}_N$ with character group $X^\bullet(T_N)= \oplus_{i=1}^N e_i$ and then form the maximal torus $mT_N$ of $m \mr{GL}_N$ making the following sequence exact:
\[
1 \to mT_N \to T_N \times \mathbf{G}_m \xrightarrow{(\det^{-1}, [m])} \mathbf{G}_m \to 1.
\]
Then $X^\bullet(mT_N)$ is the cokernel
\[
\frac{X^\bullet(T_N) \oplus \Z f}{\Z(-\sum_{i=1}^n e_i + mf)} \xrightarrow{\sim} X^\bullet(mT_N),
\]
which can alternatively be identified as 
\begin{align*}
\frac{\bigoplus_{i=1}^N \Z e_i \oplus \Z f}{\Z(-\sum_{i=1}^N e_i + mf)} &\xrightarrow{\sim} \bigoplus_{i=1}^N \Z e_i + \Z \frac{\sum_1^n e_i}{m} \\
e_i &\mapsto e_i \\
f &\mapsto \frac{\sum_i e_i}{m}.
\end{align*}
\endproof
\begin{lemma}
There is a functorial transfer from cuspidal automorphic representations of $\mr{GL}_N(\af)$ to cuspidal automorphic representations of $m\mr{GL}_N(\af)$ corresponding to the L-homomorphism $\mr{GL}_N(\CC) \to \mr{GL}_N(\CC)/\mu_m$. A mixed-parity (see Definition \ref{mixedparitydef}) automorphic representation of $\mr{GL}_N(\af)$ inflates to an L-algebraic automorphic representation of $2\mr{GL}_N(\af)$.
\end{lemma}
\proof
Cuspidal automorphic forms on $\mr{GL}_N(\af)$ can be regarded as functions on $m\mr{GL}_N(F)\backslash m\mr{GL}_N(\af)$ by inflation.\footnote{The map is neither injective nor surjective, so I don't know whether to call this inflation or restriction.} Given a $\Pi$ on $\mr{GL}_N(\af)$ we define its functorial transfer to be an irreducible constituent of the corresponding space of functions on $m\mr{GL}_N(\af)$; this inflation is $m\mr{GL}_N(\af)$-equivariant in the obvious way, so at the level of isomorphism classes of irreducible admissible representations, everywhere locally we are just restricting from $\mr{GL}_N(F_v)$ to (the image in $\mr{GL}_N(F_V)$ of) $m\mr{GL}_N(F_v)$. That this corresponds to the desired reduction modulo $\mu_m$ of L-parameters follows from \cite[Corollary 3.1.6]{stp:variationsarxiv}.

For the second claim, note that the surjection of maximal tori $2T_N \onto T_N$ (writing $T_N$ for a maximal torus of $\mr{GL}_N$) induces, in suitable coordinates, the inclusion of character groups
\[
\bigoplus_1^N \Z e_i \into \bigoplus_1^N \Z e_i + \Z \left( \frac{\sum_1^N e_i}{2} \right).
\]
\endproof
Given such a $\Pi$ on $\mr{GL}_N(\af)$, we will write $m \Pi$ for one of its automorphic transfers. To give a more explicit sense of what's going on, we give a couple of examples of the relevant local representation theory:
\begin{lemma}
Let $v$ be a finite place of $F$, and let $\Pi_v$ be an irreducible admissible smooth representation of $\mr{GL}_N(F_v)$. Let $H_v$ denote the image of $m\mr{GL}_N(F_v)$ in $\mr{GL}_N(F_v)$.
\begin{enumerate}
\item Suppose $(m, N)=1$. Then the restriction $\Pi_v|_{H_v}$ is irreducible.
\item Suppose $m=N$, in which case the natural map $\mr{SL}_N \times \mathbf{G}_m \xrightarrow{\sim} N \mr{GL}_N$ is an isomorphism. Then the irreducible constituents of $\Pi_v|_{H_v}$ are in bijection with those of $\pi_v|_{\mr{SL}_N(F_v)}$. Their number can grow arbitrarily large as $N \to \infty$.
\end{enumerate}
\end{lemma}
\proof
In all cases, $\Pi_v|_{H_v}$ is multiplicity-free, by the corresponding result for $H_v= \mr{SL}_N(F_v)$ (see \cite[Theorem 1.3]{adler-prasad}). The number of irreducible constituents of this $H_v$-representation is then (see \cite[Lemma 2.1(d)]{gelbart-knapp:l-indistinguishability}) the order of the group
\[
X_{H_v}(\Pi_v)= \{ \chi \colon \mr{GL}_N(F_v)/H_v \to \CC^\times: \Pi_v \otimes \chi \cong \Pi_v \}.
\]
Since $\mr{GL}_N(F_v)/H_v \cong F_v^\times/(F_v^\times)^m$ has exponent $m$, and (taking central characters) $\chi^N=1$, we conclude that $X_{H_v}(\Pi_v)$ is trivial whenever $(m, N)=1$.

The first assertion of Part 2 is clear, noting that the transfer from $\mr{GL}_N$ to $N\mr{GL}_N$ is associated to the more familiar morphism of dual groups $\mr{GL}_N \to \mr{PGL}_N \times \mathbf{G}_m$, the second projection simply being the determinant. The second assertion (again, just for intuition) is just a vague statement that is part of a much more precise description of L-packets for $\mr{SL}_N(F_v)$ (see \cite{gelbart-knapp:l-indistinguishability}). 
\endproof
If one were willing to work harder, these local results could be used to understand automorphic multiplicities for $m\mr{GL}_N$, similarly to the case of $\mr{SL}_N$.
\subsection{Examples} Now suppose that we have a geometric projective representation $\bar{\rho} \colon \gal{F} \to \mr{PGL}_N(\Qlb)$ arising from (the restriction to $\mr{SL}_N(\af)$ of) a mixed-parity, regular, ESD, cuspidal\footnote{Or isobaric, with the essential self-duality applying to each cuspidal constituent, since the proof of Theorem \ref{mixedgalois} works just as well.} representation of $\mr{GL}_N(\af)$. We ask whether it lifts geometrically to $\mr{GL}_N(\Qlb)/\{\pm 1\}$. Again, the interesting case will be $N$ even. In the notation of Proposition \ref{reallift}, $G= \mr{SL}_N$ and $\tG= 2\mr{GL}_N$, $Z_G= \mu_N$ and
\[
\tZ= \{(z_1, z_2) \in Z_{\mr{GL}_N} \times \mbf{G}_m: z_1^N= z_2^2\} \cong 
\begin{cases}
\mbf{G}_m \times \Z/2 & \text{if $N$ is even;}\\
\mbf{G}_m & \text{if $N$ is odd.}
\end{cases}
\]
We can check that
\[
\coker \left( X^\bullet(\tZ)_{tor} \to X^\bullet(Z_G)_{tor} \right)
\]
is isomorphic to $\mu_N/\mu_N[2]$ (the quotient by the $2$-torsion) in either case. Thus Proposition \ref{reallift} (and the proof of Theorem \ref{mixedgalois}) implies:
\begin{cor}\label{pm1liftstotreal}
Let $F$ be a totally real field, and let $\Pi= \boxplus_{i=1}^r \sigma_i$ be a regular, mixed-parity, isobaric automorphic representation of $\mr{GL}_N(\af)$, all of whose cuspidal constituents $\sigma_i$ satisfy an essential self-duality $\sigma_i \cong \sigma_i^\vee \otimes \omega$ for some Hecke character $\omega$ of $F$. For each $\ell$, fix $\iota \colon \Qlb \xrightarrow{\sim} \CC$. Then there exists a compatible system of $\ell$-adic representations
\[
\bar{\rho}_{\iota} \colon \gal{F} \to \mr{PGL}_N(\Qlb)
\]
 associated to $\Pi|_{\mr{SL}_N(\af)}$. Moreover, for each $\ell$ (and $\iota$), $\bar{\rho}_{\iota}$ lifts to a geometric representation $\gal{F} \to \mr{GL}_N(\Qlb)/\{ \pm 1 \}$. That $\Pi$ is mixed-parity also implies $\bar{\rho}_{\iota}$ has no geometric lift to $\mr{GL}_N(\Qlb)$.
\end{cor}
\begin{rmk}
As a referee kindly pointed out, using the isomorphism $\mr{SL}_N \times \mathbf{G}_m \xrightarrow{\sim} N \mr{GL}_N$, we can immediately deduce that $\Pi$ has an associated compatible system of $\ell$-adic representations valued in the groups
\[
\mr{GL}_N(\Qlb)/\mu_N \cong \left( \mr{PGL}_N \times \mathbf{G}_m \right)(\Qlb).
\]
The projective factors are given by Corollary \ref{pm1liftstotreal}, and the compatible system of Galois characters demanded by the $\mathbf{G}_m$-factor is that associated to the (type $A_0$) central character of $\Pi$.
\end{rmk}
We would like to upgrade this to the stronger statement that there exists a compatible system of $\ell$-adic representations 
\[
\rho_{\iota} \colon \gal{F} \to \mr{GL}_N(\Qlb)/\{ \pm 1\}
\]
compatible with the local parameters of $2 \Pi$. This is rather delicate: recall that $\bar{\rho}_{\iota}$ was constructed by patching representations $\rho_L^{\star}$ over CM extensions $L$ of $F$, yielding a $\gal{F}$-representation $\tilde{\rho}$ whose projectivization gave $\bar{\rho}_{\iota}$. The subtle, yet crucial, point is that the desired $\rho_{\iota}$ (in the mixed-parity case) is \textit{not} the reduction modulo $\pm 1$ of $\tilde{\rho}$ (at least for $N=2$, this can be checked unconditionally); indeed, I do not believe that $\rho_{\iota}$ can lift, at all, to $\mr{GL}_N(\Qlb)$.\footnote{Note that there are genuine obstructions to lifting through isogeny quotients such as $\mr{GL}_N \to \mr{GL}_N/\{ \pm 1 \}$; in this case, the obstructions lie in $H^2(\gal{F}, \pm 1)$.}

\bibliographystyle{amsalpha}
\bibliography{biblio.bib}

\end{document}